\documentclass{amsart}
\usepackage[dvips]{color}
\usepackage{graphicx}

\usepackage{latexsym}
\usepackage{amssymb}
\usepackage{amsmath}
\usepackage{amsthm}
\usepackage{amscd}
\theoremstyle{plain}
\newtheorem{prop}{Proposition}
\newtheorem{thm}{Theorem}
\newtheorem{lem}{Lemma}
\newtheorem{Bob}{Definition}

\newtheorem{cor}{Corollary}
\theoremstyle{remark}
\newtheorem{rem}{Remark}

\begin{document}

\title[Hausdorff Dimension of Ergodic Measures]{Hausdorff Dimension for Ergodic Measures of Interval Exchange Transformations}
\author[J.\ Chaika]{Jon Chaika}

\address{Department of Mathematics, Rice University, Houston, TX~77005, USA}

\email{jmc5124@rice.edu}

\begin{abstract}
I show that there exist minimal interval exchange transformations with
an ergodic measure whose Hausdorff dimension is arbitrarily small, even 0. I
will also show that in particular cases one can bound the Hausdorff
dimension between $\frac 1 {2r+4}$ and $\frac 1 r$  for any r greater than 1. 
\end{abstract}
\maketitle
\section{Introduction}
 This paper is concerned with determining the possible Hausdorff dimensions for ergodic measures of
interval exchange transformations (IETs). In general, almost
all IETs (with irreducible permutation) are uniquely ergodic with
respect to Lebesgue measure [V2], [M] (see also [B1] for an elementary proof) and therefore their only
ergodic measure has Hausdorff dimension 1. In fact, since a minimal
$d$ interval IET has at most [$\frac d 2 $] (probability) ergodic measures [V1], [Ka], the
smallest number of intervals for which one can non-trivially consider the Hausdorff dimension of ergodic measures is 4.

\vspace{2mm}

 In a celebrated 1977 paper Michael Keane provided a
method for constructing minimal 4 IETs that are
not uniquely ergodic [K]. This followed an earlier construction of Keynes and Newton, who showed a minimal  5 IET could be not uniquely ergodic [KN]. Using Keane's construction I will present
several results on the Hausdorff dimensions of ergodic measures of
non-uniquely ergodic minimal IET's. By the estimate in the previous paragraph, minimal 4-IETs have at most 2= [$\frac 4 2 $] ergodic measures. I will restrict my attention to when one of the ergodic measures is Lebesgue
measure and the other is singular with respect to Lebesgue
measure. If this is not the situation, then both of
the measures are absolutely continuous with respect to Lebesgue and hence have Hausdorff dimension 1. (Lebesgue measure is preserved and therefore can be written as a combination of ergodic measures.)

\vspace{1mm}

Some of
the results contained in the paper are:

\vspace{2mm}

1) That one can build minimal IETs that have an ergodic measure with arbitrarily small (even 0)
Hausdorff dimension (see Theorem 2 and Corollary 1).

\vspace{1mm}

 2) That one can place the Hausdorff dimension
of particular IETs between $\frac 1 r$ and $\frac 1 {2r+4}$ for
$r\geq 1$.

\vspace{1mm}

 3) That a large class of Keane's non-uniquely ergodic
IETs have both ergodic measures having Hausdorff dimension 1 (see Remark 1).

\vspace{2mm}

 Prior to this work it was known that if the lengths of the intervals are algebraic numbers then the Hausdorff Dimension is greater than 0 [B2]. However, it was not known if the the Hausdorff dimension could ever be less than 1. 
 
\vspace{2mm}

The plan of the paper is as follows: The second section reviews IETs and Keane's paper. The third section introduces terminology and notation that will be used throughout the remainder of the paper. The fourth section states the theorems of the paper. The fifth section presents preliminary lemmas. The sixth section presents proofs of theorems and a couple of results that they imply.

\vspace{3mm}

\section{Review of IETs and Michael Keane's Results}

\vspace{1.5mm}
\begin{Bob} Given $L=(l_1,l_2,...,l_n)$
where $l_i>0$,  $l_1+...+l_n=1$ we can obtain n subintervals of the unit
interval $I_1=[0,l_1) ,
I_2=[l_1,l_1+l_2),...,I_n=[l_1+...l_{n-1},1)$. If we are also given
 a permutation on n letters $\pi$ we obtain an n \emph{Interval Exchange Transformation}  $ T \colon [0,1) \to
 [0,1)$ which exchanges the intervals $I_i$ according to $\pi$. That is, if $x \in I_j$ then \begin{center} $T(x)= x - \underset{k<j}{\sum} l_k +\underset{\pi(k')<\pi(j)}{\sum} l_{k'}$.\end{center}
\end{Bob}
Keane's paper and this one are concerned with 4-IETs with permutation (4213).
Keane relied on the induced map on the fourth interval for his result. He showed that by choosing the lengths appropriately one could ensure that this induced map had the permutation (2431). Name these in reverse order and we once
again get a (4213) IET. Moreover, Keane showed that for any choice $ m,n \in \mathbb{N} $ one can find
 an IET whose landing pattern is given by the columns of following matrix: \\
\begin{center} $A_{m,n}=
\left(
\begin{array}{cccc}
0 & 0 & 1 &1 \\
m-1 & m & 0 & 0 \\
n & n & n-1 & n \\
1 & 1 & 1 & 1 \end{array} \right)$; $\quad$ $m, n \in \mathbb{N} =\{1, 2, ...\} $.
\end{center} For instance the second (after
renaming) subinterval of the induced map  visits the intervals of
the original IET according to the pattern $[0 \, m \, n \, 1]$ before returning
to the 4th interval. That is, it does not land in the first interval, it lands a total of $m$ times in the second interval,
 $n$ times in the third interval and returns to the fourth
interval.  One can now repeat this procedure on the 4th interval
(once again after renaming) of our induced map with a new matrix
$A_{m_2,n_2}$ and on the 4th interval of this induced map with another
matrix $A_{m_3,n_3}$ and so on. 

\vspace{2mm}

The IETs that have this property
for the matrices $A_{m_1,n_1},.., A_{m_k,n_k}$ are those contained
in the image of the 3-simplex under the map
$\bar{A}_{m_1,n_1}\bar{A}_{m_2,n_2}...\bar{A}_{m_k,n_k}$ (see Definition 7). Michael Keane showed that if we choose
our $m_i , n_i$ appropriately we get an IET with
two ergodic measures. In particular, if one chooses our lengths
according to the vector $\underset{k \to
\infty}{\lim}\bar{A}_{m_1,n_1}\bar{A}_{m_2,n_2}...\bar{A}_{m_k,n_k}e_3$ (where $e_3= [0 \, 0\, 1 \, 0]^T)$ one gets
that one of our measures is Lebesgue measure and the other is
singular with respect to Lebesgue measure. The singular measure
assigns weights according to the vector $\underset{k \to
\infty}{\lim}\bar{A}_{m_1,n_1}\bar{A}_{m_2,n_2}...\bar{A}_{m_k,n_k}e_2$ (where $e_2=[0 \, 1\, \, 0\, 0]^T$). Following
Keane we will denote the Lebesgue measure as $\lambda_3$ and the
singular measure $\lambda_2$. The conditions Keane gives are as
follows (for notation see Definition 7):
\begin{thm} If one chooses $3(n_k +1)\leq m_k \leq \frac 1 2 (n_{k+1}+1)$
and $n_1\geq 10$ then the IET with lengths determined by the vector $\underset{k \to \infty}{\lim} \bar{A}_{m_1,n_1}\bar{A}_{m_2,n_2}...\bar{A}_{m_k,n_k}e_3$ is not
uniquely ergodic.
\end{thm}
 This is Theorem 5 of [K] and these conditions are assumed to be satisfied for the remainder of the paper. The following two
results of Keane are key in proving this result and will be used
in this paper (for notation see Definition 2):
\begin{lem} $\frac{\lambda_3(I_2^{(k)})}{\lambda_3(I^{(k)})} \leq \frac {2m_k} {(n_{k+1}+1)(n_k +1)}$.
\end{lem}
This result is in the proof of Lemma 3 of [K].
\begin{lem} $\frac{\lambda_2(I_2^{(k)})}{\lambda_2(I^{(k)})} \geq \frac 1 3 $.
\end{lem}
This is Lemma 4 of [K].

\vspace{3mm}

\section{Definitions and Notation}

\vspace{1.5mm}

\begin{Bob} $I^{(k)}$ denotes the k$^{\text{th}}$ induced interval
and $I_j ^{(k)}$ denotes the $j^{\text{th}}$ subinterval of $I^{(k)}$.
\end{Bob}
Note: $I_4^{(k)}=I^{(k+1)}$. Also note that in [K]
this notation is reversed, so in his paper $I_k ^{(j)}$ is the $j^{\text{th}}$
subinterval of the $k^{\text{th}}$ induced interval.

\vspace{1mm}

\begin{Bob} $|v|_{1}$ is the sum of the absolute values of the
entries of the vector $v$.
\end{Bob}

\vspace{1mm}

\begin{Bob} $B_k=A_{m_1,n_1}A_{m_2,n_2}...A_{m_k,n_k}$.
\end{Bob}
This matrix describes the travel of the subintervals of
$I^{(k)}$ until they land in $I^{(k)}$ again. That is, the number of times each subinterval of $I^{(k)}$ lands in our initial subintervals before returning to $I^{(k)}$.
\vspace{1mm}
\begin{Bob} $b_{t,i}$ denotes $|B_te_i|_{1}$.
\end{Bob}
As above, $e_i$ denotes the column vector where the $i^{\text{th}}$ entry is 1 and all other entries are 0.
\vspace{1mm}

\begin{Bob} $O(I_j^{(k)})$ denotes the union of images of $I_j^{(k)}$ that
$B_ke_j$ counts.
\end{Bob}
That is, $O(I_j^{(k)})= \underset{l=0}{\overset{{b_{k,j}-1}}{\cup}} T^l(I_j^{(k)})$.

\vspace{1mm}

\begin{Bob} Let $\bar{A}_{m_k , n_k}(v)=\frac {A_{m_k,n_k}(v)}
{|A_{m_k,n_k}(v)|_{1}}$.
\end{Bob}
$\bar{A}_{m_k , n_k}$ maps vectors in the unit 3-simplex to vectors in the unit 3-simplex. This ensures that the measures obtained by the conditions of the theorems are probability measures. 

Note: $\bar{A}_{m_k , n_k}(u+v) \neq \bar{A}_{m_k , n_k}(u)+\bar{A}_{m_k , n_k}(v)$ in general.
\vspace{1mm}

\begin{Bob} $S=\underset{k=1}{\overset{\infty}{\cap}} \underset{r=k}{\overset{\infty}{\cup}} O(I_2^{(r)})$.
\end{Bob}
This is the set of all points which lie in 
$O(I_2^{(k)})$ for infinitely many $k$.

\vspace{1mm}

\begin{Bob} If $M \subset [0,1)$ is set, $H_{\text{dim}}(M)$ denotes the \emph{Hausdorff dimension} of~$M$. \end{Bob}
For a definition of Hausdorff dimension and an introduction to it
see [F].
\begin{Bob} The \emph{Hausdorff dimension of a probability measure
$\mu$} is \begin{center}$H_{\text{dim}}(\mu)=
\inf\{H_{\text{dim}}(M) \colon M \text{ is Borel and } \mu(M)=1 \}$.\end{center}
\end{Bob}

\vspace{3mm}

\section{New Results}

\vspace{1.5mm}

\begin{thm} If an IET has lengths determined by the vector
\begin{center} $\underset{k \to \infty}{\lim} \bar{A}_{m_1,n_1}\bar{A}_{m_2,n_2}...\bar{A}_{m_k,n_k}e_3$ \end{center} and there exists N such that $n_{k+1}\geq (b_{k,2})^r 2^{rk} m_k$ for all $k \geq N$, then the Hausdorff dimension of $\lambda_2$, the other ergodic measure,  is less than or equal to $\frac 1 r$.
\end{thm}
The condition for Theorem 2  along with Lemma 1 implies that \\ $\lambda_3
(I_2^{(k)}) \leq \frac 1 {(b_{k,2})^r 2^{rk} }$ for $k \geq N$. This fact is crucial for the proof.
 \begin{thm}  If an IET has lengths determined by the vector
\begin{center} $\underset{k \to \infty}{\lim} \bar{A}_{m_1,n_1}\bar{A}_{m_2,n_2}...\bar{A}_{m_k,n_k}e_3$ \end{center} and there exists N such that $b_{k+1,2}\leq (b_{k,2})^r$, $m_k\geq k^2 n_k$  for all $k\geq N$, then the Hausdorff dimension of $\lambda_2$, the other ergodic measure, is greater than or equal to~$\frac 1 {2r}$.

\end{thm}
As the next theorem suggests, there is a gap between the
$r$ in Theorem 2 and in Theorem 3. In general one can have
$r_{2} + 2 \geq r_{3}$ (where $r_2$ is the $r$ in Theorem 2 and $r_3$ is the $r$ in Theorem 3). This is done by setting $n_k=k^2(b_{k-1,2})^r2^{r(k-1)}m_{k-1}$ and $m_k=k^2n_k$. Theorem 2 provides the upper bound and Theorem 3 provides the lower bound. In particular,

\begin{thm} If an IET has lengths determined by the vector
\begin{center} $\underset{k \to \infty}{\lim} \bar{A}_{m_1,n_1}\bar{A}_{m_2,n_2}...\bar{A}_{m_k,n_k}e_3$, \end{center} with $n_k=9^{4^{k-1}}$ and $m_k=9^{4^{k-1}+k}$,
then $\frac 1 8 \leq H_{\text{dim}} (\lambda_2) \leq \frac 1 2 $.
\end{thm}

\vspace{3mm}

\section{Preliminary Lemmas}

\vspace{1.5mm}

First, a strengthening of Lemma 2.
\begin{lem} $\lambda_2(O(I_2^{(k)}))=b_{k,2} \lambda_2(I_2 ^{(k)})$ is
greater than $\frac 1 3 $ for all $k \geq 0$.
\end{lem}
\noindent Proof: I begin by showing that $b_{k,2} \geq b_{k,i}$ by comparing the entries of $B_ke_2$ and $B_ke_i$. $b_{k,2}>b_{k,1}$ because the
 second entry of $A_{m_k, n_k} e_2=m_k > m_k-1$ and $m_k-1$ is the second entry of $A_{m_k, n_k} e_1$.
$A_{m_k,n_k}e_2$ agrees with
$A_{m_k,n_k}e_1$ in all other entries. $b_{k,2} \geq b_{k,j}$ for $j=3,4$ because $A_{m_k,n_k}e_2 \geq A_{m_k,n_k} e_j$ in
all entries but the first and
$m_kA_{m_{k-1},n_{k-1}}e_2>A_{m_{k-1},n_{k-1}}e_1$ in all entries (the second entry of $A_{m_k,n_k}e_j$ is 0 and the second entry
of $A_{m_k,n_k}e_2$ is $m_ke_2$ and also the first entry of $A_{m_k,n_k}e_j =1$). This argument shows that $A_{m_{k-1},n_{k-1}}A_{m_k,n_k}e_2$ has each entry greater than or equal to the corresponding entries of $A_{m_{k-1},n_{k-1}}A_{m_k,n_k}e_j$ for $j=3,4$.

We also have that $\lambda_2(I_2^{(k)}) > \frac 1 3
\lambda_2(I^{(k)})$ by Lemma 2. Therefore, because our IET is minimal,
we have,  
\begin{align*}
 1=\lambda_2([0,1])&=\lambda_2(O(I_1 ^{(k)}))+\lambda_2(O(I_2
^{(k)}))+\lambda_2(O(I_3 ^{(k)}))+\lambda_2(O(I_4
^{(k)}))\\ &=b_{k,1}\lambda_2(I_1^{(k)}) +
b_{k,2}\lambda_2(I_2^{(k)})+ b_{k,3}\lambda_2(I_3^{(k)}) +
b_{k,4}\lambda_2(I_4^{(k)})  ,
\end{align*} 
and so $\lambda_2(O(I_2 ^{(k)}))=\frac 1
{3}$. In fact, $\lambda_2(O(I_2^{(k)}))\geq \frac {\lambda_2(I_2
^{(k)})}{\lambda_2(I^{(k)})}$.

\vspace{3mm}

An immediate consequence of this
lemma is that $\lambda_2(S) \geq \frac 1 3 $. In fact,

\vspace{2mm}

\begin{lem} $\lambda_2(S)=1$.
\end{lem}
\noindent Proof: By ergodicity of $\lambda_2$, it suffices to show that $\lambda_2(\underset{r=1}{\overset{\infty}{\cap}}T^r(S))\geq \frac 1 3 $. Observe that $T(O(I_2^{(k)}))\cap O(I_2^{(k)})$ is at most missing the last image of $I_2^{(k)}$ contained in $T(O(I_2^{(k)}))$. So $\underset{k \to \infty}{\lim} \lambda_2(O(I_2^{(k)}) \cap  T^d (O(I_2^{(k)})))= \underset{k \to \infty}{\lim} \lambda_2( O(I_2^{(k)}))$ and therefore:
\begin{center} $\lambda_2(\underset{r=1}{\overset{\infty}{\cap}}T^r(S))= \underset{d \to \infty} {\lim} \lambda_2 (\underset{r=1}{\overset{d}{\cap}}T^r(S)) \geq \frac 1 3 $. \end{center}

\begin{lem} Under the conditions of Theorem 3 $ \quad \frac{\lambda_2(I_2 ^{(k)})}{\lambda_2(I^{(k)})}>
\frac{k^2}{k^2+4}$.
\end{lem}
This lemma follows by induction.

\vspace{2mm}

An immediate consequence of this lemma and the proof of
Lemma 3 is that $\lambda_2 (O(I_2^{(k)}))\geq \frac {k^2}{k^2+4}$
too.

It follows from this and Borel-Cantelli that $\lambda_2$ almost every point
is in $O(I_2^{(k)})$ for all $k$ large enough. This leads to the
following lemma:

\vspace{2mm}

\begin{lem} Under the conditions of Theorem 3, for $\lambda_2$ almost all points $x$, \mbox{$T^{-r}(x) \in I_2
^{(k)}$} and $T^s(x) \in I_2 ^{(k)}$ for some $0\leq r,s\leq b_{k,2}$ for all
but finitely many $k$.
\end{lem}

\vspace{2mm}

 This lemma says that for $ \lambda_2$ almost every $x$ there exists $N_x$
such that \mbox{$T^i (x) \in O(I_2^{(k)})$} for $k \geq N_x $ and $0 \leq i
\leq b_{k,2}$.

\vspace{2mm}

\begin{lem} $\lambda_3(I_1^{(k)})\geq \lambda_3 (I^{(k+1)}) \frac{ \lambda_3
(I_3^{(k+1)})}{\lambda_3(I^{(k+1)} ) }$.
\end{lem}

\noindent Proof: Observe that $\lambda_3 (I_1^{(k)})=
\lambda_3(I^{(k)})\frac{\lambda_3(I_3^{(k+1)} \cup I_4^{(k+1)})}
{\lambda_3(I^{(k)})}$. Also, $b_{k,2} \geq b_{k,i}$ (see Lemma
3) and $T$ is minimal implies $\lambda_3 (I^{(k)}) \geq \frac 1 {b_{k,2}}$.

\vspace{2mm}

This lemma is similar to Lemma 1 in [K].

\vspace{2mm}

\begin{lem} The images of $I_2^{(k)}$ in $O(I_2^{(k)})$ are never immediately
adjacent.
\end{lem}
\noindent Proof:  By Keane's construction $I_2^{(k)}$ will always be
bordered on both sides by $I_1^{(k)}$ or on one side by $I_1^{(k)}$
and the other by $I_4^{(k)}$ or on one side by $I_1^{(k)}$ and the
other by $I_3^{(k)}$. This is because the image of $I_2^{(k)}$ that
are inside a subinterval of $I^{(k-1)}$ have this property. Also the
only subinterval of $I^{(k-1)}$ which has an image of $I_2^{(k)}$ on
its boundary is $I_2^{(k-1)}$ (its left and right boundary are both
images of $I_2^{(k)}$). The result follows by induction on k. Just
for reference, $I_1^{(k-1)}$'s boundary blocks are $I_3^{(k)}$ and
$I_4^{(k)}$, $I_3^{(k-1)}$'s are $I_1^{(k)}$ and $I_4^{(k)}$ and
$I_4^{(k-1)}$'s boundary blocks are $I_4^{(k)}$ and $I_1^{(k)}$.

\vspace{3mm}

\section{Proofs of the Theorems}

\vspace{1.5mm}

\begin{prop} Under the conditions of Theorem 3, $\lambda_2$ a.e.
x satisfies \begin{center}$\underset{s \to \infty} {\lim \inf} \, s^{2r+\epsilon}
|T^s(x)-x|=\infty$.\end{center}
\end{prop}

\noindent Proof: Given $x$, pick $k'$ so that we have $x$ satisfying
Lemma 6 for all $k>k'$. By Lemma 8 this means that for $s \leq
b_{k,2}$ we have $|T^s (x)-x|\geq \underset{i \neq 2} {\min} \, \{
\lambda_3(I_i^{(k)}) \} $. This is because all of the images of $x$
lie in separate images of $I_2^{(k)}$, which are separated by the
image of some $I_i^{(k)}$ for $i \neq 2$. $I_1^{(k)}$ has the
smallest $\lambda_3$ measure of these subintervals (subintervals that are not
$I_2^{(k)}$). $\lambda_3(I_1^{(k)}) $ gives a lower bound. By
Lemma 7, the fact that $b_{k,2} \geq b_{k,i}$ and the fact that
$\underset {k \to \infty}{\lim} \frac {\lambda_3 (I_3 ^{(k)})
}{\lambda_3 (I^{(k)})}=1$, it follows that $\lambda_3 (I_1^{(k)}) \geq \frac
1 {(b_{k+1,2})^{1+\epsilon}}$, eventually. (Indeed, $b_{k,2} \geq b_{k,i}$ so $\lambda_3(I^{(k)})\geq \frac 1 {b_{k,2}}$.) Thus, $|T^s (x)-x|\geq
\frac 1 { (b_{k+1,2})^{1+\epsilon}}$. So if the conditions of
Theorem 3 are satisfied and 
$b_{k,2}\leq s \leq b_{k+1,2}$, then we have
\begin{align*} s^{2r+4r\epsilon}|T^s(x)-x| &\geq
(b_{k,2})^{2r+4r\epsilon}  \lambda_3(I_1^{(k+1)})\geq \\  
(b_{k,2})^{2r+4r\epsilon} \frac 1 { (b_{k+2,2})^{1+\epsilon}} & \geq
(b_{k,2})^{2r+4r\epsilon} \frac 1 { ((b_{k,2})^{2r}) ^{1+\epsilon }
}=\\  (b_{k,2})^{2r+4r \epsilon} \frac 1 {(b_{k,2})^{2r+2r \epsilon
}} & =(b_{k,2})^{2r\epsilon}, 
 \end{align*}  which goes to infinity with $k$.

Theorem 3 can now be proved with the assistance of Theorem 1.3 in
 [B2]. Put in the
language of this paper it states:
\begin{thm} If the Hausdorff dimension of an
invariant measure $\mu$ for a dynamical system $T:[0,1) \to [0,1)$
is less than  $\alpha$ then $\lim \inf \{n^{\frac 1
{\alpha}}|T^n(x)-x|\}=0$ for $\mu$ almost every $x$.
\end{thm}
 This proves Theorem 3 because it shows
that $H_{\text{dim}}(\lambda_2) \geq \frac {1} {2r+\epsilon}$ for
any $\epsilon$.

\begin{rem}
Theorem 3 also shows that if $m_k \geq k^2 n_k$ then, unless
one stipulates much faster growth than Keane does,
$H_{\text{dim}}(\lambda_2)=1$. This is because if $n_k$ grows
exponentially so does $m_k$. This implies that $b_{k,2}$ grows like $c^{k^2}$, which means that for any $\epsilon>0$ eventually
$(b_{k,2})^{1+\epsilon}>b_{k+1,2}$.
\end{rem}

\vspace{3mm}

\noindent Proof of Theorem 2:  By Lemma 4 it suffices to show $H_{\text{dim}}(S)
\leq \frac 1 r $. Observe that covering $O(I_2^{(k)})$ with images of $I_2^{(k)}$, performing
the Hausdorff $\frac 1 r $ dimensional estimate gives a number
less than $b_{k,2}(\frac 1 {(b_{k,2})^r 2^{rk}})^{\frac 1 r }=2^{-k}$.
By summing $\lambda_2(O(I_2^{(k)}))$ over $k\geq L$ the Hausdorff
$\frac 1 r $ dimensional measure of $S$ is less than $2^{-L+1}$ for
any $L$. So $H_{\text{dim}}(S)\leq \frac 1 r$.

\begin{cor} If the conditions of Theorem 2 are satisfied and additionally, \mbox{$n_{k+1}=(b_{k,2})^k$}, then 
$H_{\text{dim}}(\lambda_2)=0$.
\end{cor}

\noindent Proof: The conditions of Theorem 2 are satisfied for any $r$ by
picking $N$ big enough.

\vspace{3mm}

In order to prove Theorem 4, we obtain coarse estimates on $b_{k,2}$.
\begin{lem} $b_{k,2}\geq m_1m_2...m_k$.
\end{lem}
 Indeed, the second entry of $B_ke_2$ is bigger than this.
\begin{lem} $b_{k,2}\leq (m_1+n_1+1)(m_2+n_2+1)...(m_k+n_k+1)$.
\end{lem}
\noindent Proof: Observe that $A_{m,n}e_2 \geq A_{m,n}e_j$, which means $|A_{m,n}v|_{1}\leq
|v|_{1}|A_{m,n}e_2|_{1}$. The lemma follows from this fact and
induction.
\\ \\
Proof of Theorem 4: Observe $(m_i+n_i+1)\leq 2m_i$. So by Lemmas 9 and 10 under the
conditions of Theorem 4 we have \begin{center} $9^{(4^{k}-1)\frac 1 3 +\frac
{k(k+1)} {2} } \leq b_{k,2} \leq 2^k 9^{(4^{k}-1)\frac 1 3 +\frac
{k(k+1)} {2} }$.\end{center} In view of the fact that \begin{center}$(2^k 9^{(4^{k}-1)\frac 1 3 +\frac
{k(k+1)} {2} })^2 2^{2k} 9^{4^{k-1}+k}<9^{4^{k}}$ \end{center} for large $k$, the conditions of Theorem 2 are satisfied
with $r_1=2$. Also notice that
\begin{center}$(9^{(4^{k}-1)\frac 1 3 +\frac {k(k+1)} {2} })^{4+\epsilon}\geq
2^{k+1}9^{(4^{k+1}-1)\frac 1 3 +\frac {(k+1)(k+2)} {2} }$ \end{center}
for large $k$, satisfying the conditions for Theorem 3 with
$r_2=4+\epsilon$. This gives Theorem 4.

\section{Acknowledgments}
I would like to thank my advisor, Professor M. Boshernitzan, for
posing this problem to me and his numerous helpful conversations. I
would like to thank Professor A. Bufetov for suggestions of results
to include in this paper. I would like to thank Professor C.
Ulcigrai for pointing out a gap in the initial proof. I would also like to thank Professor R. Grigorchuk, Professor Y. Vorobets and the anonymous referee for suggestions to improve the clarity of the paper.

\vspace{3mm}

\section{References}
\noindent [B1] Boshernitzan, M.: A Condition for minimal interval exchange maps to be uniquely ergodic. Duke J. of Math. 52 (3) (1985) 723-752.

\noindent [B2] Boshernitzan, M.: Quantitative recurrence results. Invent.
Math. 113 (3)(1993) 617-631.

\noindent [F] Falconer, K: Fractal Geometry. Mathematical
Foundations and Applications. Wiley, Chichester, 1990.

\noindent [Ka] Katok, A: Invariant measures of flows on oriented surfaces. Sov. Math Dokl. 14 (1973) 1104-1108.

\noindent [K] Keane, M: Non-ergodic interval exchange
transformations, Israel J. Math. 26 (2) (1977)
188-196.

\noindent [KN] Keynes, H. B., Newton, D. : A ``minimal'', non-uniquely ergodic interval exchange transformation. Math. Z. 148 (1976) 101-105.

\noindent [M] Masur, H: Interval exchange transformations and
measured foliations. Ann. of Math. (2) 115 (1982) 168-200.

\noindent [V1] Veech, W: Interval exchange transformations. J.
D'Analyse Math. 33 (1978) 222-272.

\noindent [V2] Veech, W: Gauss measures for transformations on the space of interval exchange maps. Ann. of Math. (2) 115 (1982) 201-242.

\end{document}